\documentclass[11pt]{amsart}
\usepackage[all]{xy}

\usepackage{amssymb, amsmath}
\usepackage{eufrak}
\usepackage{amsfonts}

\newcommand{\F}{F}

\newcommand{\X}{\mathcal{X}}
\newcommand{\Y}{\mathcal{Y}}

\newcommand{\C}{\mathbf{C}}
\newcommand{\Z}{\mathbf{Z}}

\newcommand{\Q}{\mathbf{Q}}

\newcommand\jeden {1\hskip-3.5pt1}

\def\qed{\hfill $\Box$}
\def\proof{\noindent {\sl Proof} :\;  }

\newtheorem{theorem}{Theorem}

\newtheorem{lemma}{Lemma}

\newtheorem{remark}{Remark}

\begin{document}

\title[Note on Chern classes]
{ \bf 
A note on Chern-Schwartz-MacPherson class}
\author[T.~Ohmoto]{Toru Ohmoto}

\address[T.~Ohmoto]{Department of Mathematics,
Faculty of Science,  Hokkaido University,
Sapporo 060-0810, Japan}
\email{ohmoto@math.sci.hokudai.ac.jp}
%
%
\maketitle
\begin{abstract}
This is a note about 
the Chern-Schwartz-MacPherson class for certain algebraic stacks 
which has been introduced in \cite{O1}. 
We also discuss other singular Riemann-Roch type formulas in the same manner. 
\end{abstract}

\section{Introduction}

In this note we state a bit detailed account about 
MacPherson's Chern class transformation $C_*$ for quotient stacks defined in  \cite{O1}, 
although all the instructions have already been made in that paper. 
Our approach is also applicable for 
other additive characteristic classes, e.g., 
 Baum-Fulton-MacPherson's Todd class transformation \cite{BFM} 
 (see \cite{EG2, BZ} for the equivariant version) and more generally 
 Brasselet-Sch\"urmann-Yokura's Hirzebruch class transformation \cite{BSY} (see section 4 below). 
Throughout we work over the complex number field $\C$ 
or a base field $k$ of characteristic $0$. 

We begin with recalling $C_*$ for  schemes and algebraic spaces. 
These are spaces having trivial stabilizer groups. 
In following sections we will deal with quotient stacks having affine stabilizers, 
in particular,  `(quasi-)projective' Deligne-Mumford stacks 
in the sense of Kresch \cite{Kresch}.

\subsection{Schemes}
For the category of quasi-projective schemes $U$ and proper morphisms, 
there is a unique natural transformation from the constructible function functor to 
the Chow group functor,  $C_*: \F(U) \to A_*(U)$,
 so that it satisfies the normalization property:  
 $$C_*(\jeden_U)=c(TU) \frown [U] \in A_*(U) \quad \mbox{ if $U$ is smooth. }$$
This is called the {\it Chern-MacPherson transformation}, 
see MacPherson \cite{Mac} in complex case ($k=\C$) and Kennedy  \cite{Ken} 
in more general context of $ch(k)=0$. 
Here the naturality means the commutativity  $f_*C_*=C_*f_*$ of $C_*$ 
with pushforward of proper morphisms $f$. 
In particular,  for proper $pt: U \to pt (={\rm Spec}(k))$, 
the ($0$-th) degree of $C_*(\jeden_U)$ is equal to 
the Euler characteristic of $U$: $pt_*C_*(\jeden_U)=\chi(U)$ 
(as for the definition of $\chi(U)$ in algebraic context, see \cite{Ken, Joyce}). 

As a historical comment, 
Schwartz \cite{Schwartz} firstly studied 
a generalization of the Poincar\'e-Hopf theorem for 
complex analytic singular varieties by introducing   
a topological obstruction class for certain stratified vector frames, 
which in turn coincides with MacPherson's Chern class \cite{BS}. 
Therefore, $C_*(U):=C_*(\jeden_U)$  is  usually 
called the {\it Chern-Schwartz-MacPherson class} (CSM class) of a possibly singular variety $U$.

To grasp quickly what the CSM class is, there is a convenient way  
due to Aluffi \cite{Aluffi1, Aluffi2}. 
Let $U$ be a singular variety and 
$\iota: U_0 \hookrightarrow U$ a smooth open dense reduced subscheme.  
By means of resolution of singularities, 
we have a birational morphism $p: W \to U$ so that $W=\overline{U_0}$ is smooth 
and $D=W-U_0$ is a  divisor with smooth irreducible components $D_1, \cdots , D_r$ 
having normal crossings. 
Then by induction on $r$ and properties of $C_*$ it is shown that 
$$C_*(\jeden_{U_0}) =p_*\left(\frac{c(TW)}{\prod (1+D_i)} \frown [W]\right) \in A_*(U).$$
(Here $c(TW)/\prod (1+D_i)$ is equal to the total Chern class of dual to 
$\Omega_W^1(\log D)$ of differential forms with logarithmic poles along $D$). 
By taking a stratification $U=\coprod_{j} U_j$, we have 
$C_*(U) = \sum_j C_*(\jeden_{U_j})$. 
Conversely, 
we may regard this formula as an alternative definition of CSM class, 
see \cite{Aluffi1}.


\subsection{Algebraic spaces} 
We extend $C_*$ to the category of arbitrary schemes or algebraic spaces 
(separated and of finite type). 
To do this, we may generalize Aluffi's approach, or 
we may trace the same inductive proof by means of Chow envelopes 
(cf. \cite{Kimura}) 
of the singular Riemann-Roch theorem for arbitrary schemes \cite{FG}. 

Here is a short remark. 
An {\it algebraic space} $X$  is 
a stack over $Sch/k$, under \'etale topology, whose stabilizer groups are trivial: 
Precisely, there exists a scheme $U$  (called an {\it atlas}) 
and a morphism of stacks $u: U \to X$ 
such that for any scheme $W$ and any morphism $W \to X$ 
the (sheaf) fiber product $U\times_X W$ exists as a scheme, 
and the map $U\times_X W \to W$ is an \'etale surjective morphism of schemes. 
In addition, $\delta: R:=U\times_X U \to U \times_k U$ is quasi-compact, 
called the {\it \'etale equivalent relation}. 
Denote by $g_i: R \to U$ (i=1,2) the projection to each factor of $\delta$. 
The Chow group $A_*(X)$  
 is defined using an \'etale atlas $U$
  (Section 6 in \cite{EG}). 
In particular, letting $g_{12*}:=g_{1*}-g_{2*}$,  
$$\xymatrix{
A_*(R) \ar[r]^{g_{12*}} & A_*(U)  \ar[r]^{u_*} & A_*(X) \ar[r] & 0 
}$$
is exact (Kimura \cite{Kimura}, Theorem 1.8). 
Then the CSM class of $X$ is given by 
$C_*(X) = u_*C_*(U)$:  
In fact, if $U' \to X$ is another atlas for $X$ with the relation $R'$, 
we take the third $U''=U\times_X U'$ with $R''=R\times_X R'$,  
where $p: U'' \to U$ and $q: U''\to U'$ are \'etale and  finite. 
Chow groups of atlases modulo ${\rm Im}\, (g_{12*})$ are mutually identified 
through the pullback $p^*$ and $q^*$, and particularly, 
$p^*C_*(U) = C_*(U'')=q^*C_*(U')$, that is checked by using 
  resolution of singularities or  
the Verdier-Riemann-Roch \cite{Yokura} for $p$ and $q$. 
Finally we put $C_*: \F(X) \to A_*(X)$ by sending $\jeden_W \mapsto \iota_*C_*(W)$ 
for integral algebraic subspaces $W \stackrel{\iota}{\hookrightarrow} X$ and extending it linearly, and 
the naturality for proper morphisms 
is proved again using atlases. 
This is somewhat a prototype of $C_*$ for quotient stacks described below.

\section{Chern class for quotient stacks} 

\subsection{Quotient stacks}  
Let $G$ be a linear algebraic group acting on a scheme or algebraic space $X$. 
If the $G$-action is set-theoretically free, 
i.e., stabilizer groups are trivial, 
then the quotient $X \to X/G$ always exists as a morphism of algebraic spaces 
 (Proposition 22, \cite{EG}). 
Otherwise, in general we need the notion of quotient stack. 

The {\it quotient stack} $\X=[X/G]$ is 
a (possibly non-separated) Artin stack over $Sch/k$, under fppf topology 
(see, e.g., Vistoli \cite{Vistoli}, G\'omez \cite{Gomez} for the detail): 
An object of $\X$ is a family of $G$-orbits in $X$ 
parametrized by a scheme or algebraic space $B$, 
that is, a diagram $B \stackrel{q}{\leftarrow} P \stackrel{p}{\rightarrow} X$ 
where  $P$ is an algebraic space, 
$q$ is a $G$-principal bundle and $p$ is a $G$-equivariant morphism. 
A morphism of $\X$ is a $G$-bundle morphism $\phi: P \to P'$ so that $p'\circ \phi = p$, 
where $B' \stackrel{q'}{\leftarrow} P' \stackrel{p'}{\rightarrow} X$ is another object. 
Note that there are possibly many non-trivial automorphisms $P \to P$ over 
the identity morphism $id: B \to B$, 
which form the stabilizer group associated to the object 
(e.g.,  the stabilizer group of a `point' ($B=pt$) is non-trivial in general). 
A morphism of stacks $B \to \X$ naturally corresponds to an object $B \leftarrow P \to X$, 
that follows from Yoneda lemma: 
In particular there is a morphism (called {\it atlas}) $u: X \to \X$ 
corresponding to the diagram $X \stackrel{q}{\leftarrow} G\times X \stackrel{p}{\rightarrow} X$, 
being $q$ the projection to the second factor and $p$ the group action. 
The atlas $u$ recovers any object of $\X$ 
by taking fiber products: $B \leftarrow P=B \times_\X X \rightarrow X$. 

Let $f: \X \to \Y$ be a {\it proper} and  
{\it representable} morphism of quotient stacks, i.e., 
for any scheme or algebraic space  $W$ and morphism $W \to \Y$, the base change 
$\X \times_\Y W \to W$ is a proper  morphism of algebraic spaces. 
Take presentations $\X=[X/G]$, $\Y=[Y/H]$, 
and the atlases $u: X \to \X$, $u': Y \to \Y$.  
There are two aspects of $f$: \\

\noindent 
(Equivariant morphism): 
Put $B:=\X\times_\Y Y$, which naturally has a $H$-action 
so that $[B/H]=[X/G]$, $v: B \to \X$ is a new atlas, 
and $\bar{f}:B \to Y$ is $H$-equivariant: 
\begin{equation}\label{d1}
\xymatrix{
B \ar[r]^{\bar{f}} \ar[d]_{v}  & Y \ar[d]^{u'} \\
\X \ar[r]_f & \Y
}
\end{equation}
(Change of presentations): 
Let $P:=X \times_\X B$,  then the following diagram is considered as 
a family of $G$-orbits in $X$ and simultaneously  as a family of $H$-orbits in $B$, i.e., 
$p: P \to X$ is a $H$-principal bundle and $G$-equivariant, 
$q: P \to B$ is a $G$-principal bundle and $H$-equivariant: 
\begin{equation}\label{d2}
\xymatrix{
P \ar[r]^{q} \ar[d]_{p} & B  \ar[d]^{v} \\
X  \ar[r]_{u}  & \; \X.
}
\end{equation}

A simple example of such $f$  is given by 
proper $\varphi: X \to Y$ with an injective homomorphism 
$G \to H$ so that $\varphi(g.x)=g.\varphi(x)$ and $H/G$ is proper. In this case, 
$P=H \times_k X$ and $B=H\times_G X$ with 
$p: P \to X$ the projection to the second factor, $q: P \to B$ the quotient morphism.

\subsection{Chow group and pushforward} 
For schemes or algebraic spaces $X$ (separated, of finite type) with $G$-action, 
the {\it $G$-equivariant Chow gourp} $A_*^G(X)$ has been introduced  
in Edidin-Graham \cite{EG},  
and the {\it $G$-equivariant constructible function} $\F^G(X)$ in \cite{O1}. 
They are based on Totaro's algebraic Borel construction: 
Take a Zariski open subset $U$ in  
an $\ell$-dimensional linear representation $V$ of $G$ 
so that $G$ acts on $U$ freely. 
The quotient exists as an algebraic space, denoted by $U_G=U/G$. 
Also $G$ acts $X \times U$ freely, 
hence the mixed quotient $X \times G \to X_G:=X\times_G U$ 
exists as an algebraic space. 
Note that $X_G \to U_G$ is a fiber bundle with fiber $X$ and group $G$. 
Define   
 $A_n^G(X):=A_{n+\ell-g}(X_G)$ ($g=\dim G$) and $F^G(X):=\F(X_G)$ 
for $\ell \gg 0$.  
Precisely saying, we take the direct limit over all linear representations of $G$, 
see \cite{EG, O1} for the detail.

$A_n^G(X)$ is trivial for $n > \dim X$ but 
it may be non-trivial for negative $n$. Also note that 
the group $\F^G_{inv}(X)$ of {\it $G$-invariant} functions over $X$ is a subgroup of $F^G(X)$. 

Let us recall the proof that these groups are actually invariants of quotient stacks $\X$. 
Look at the diagram (\ref{d2}) above. 
Let $g=\dim G$ and $h=\dim H$. Note that $G\times H$ acts on $P$. 
Take  open subsets $U_1$ and $U_2$ of 
representations of $G$ and $H$, respectively ($\ell_i=\dim U_i\; i=1,2$) 
so that $G$ and $H$ act on $U_1$ and $U_2$ freely respectively.  
Put  $U=U_1 \oplus U_2$, on which $G \times H$ acts freely. 
We denote the mixed quotients for spaces arising in the diagram (\ref{d2})  by 
$P_{G\times H}:=P\times_{G\times H}U$, $X_G:=X\times_G U_1$ and $B_H:=B\times_H U_2$. 
Then the projection $p$ induces the fiber bundle 
$\bar{p}: P_{G\times H} \to X_G$ with fiber $U_2$ and group $H$,  
and $q$ induces $\bar{q}: P_{G\times H} \to B_{H}$ with fiber $U_1$ and group $G$. 
Thus, the pullback $\bar{p}^*$ and $\bar{q}^*$ for Chow groups  are isomorphic, 
$A_{n+\ell_1}(X_G) \simeq A_{n+\ell_1+\ell_2}(P_{G\times H}) \simeq A_{n+\ell_2}(B_H)$. 
Taking the limit, we have  the {\it canonical identification} 
$$\xymatrix{
A^G_{n+g}(X) \ar[r]^{p^*\;\; }_{\simeq\;\;} & A^{G\times H}_{n+g+h}(P)  & \ar[l]^{\;\;\simeq}_{\;\;q^*}  A^H_{n+h}(B)
}$$
(Proposition 16 in \cite{EG}). 
Note that $(q^*)^{-1}\circ p^*$ shifts the dimension by $h-g$. 
Also for constructible functions,  put the pullback $p^*\alpha:=\alpha \circ p$, 
then we have $\F^G(X) \simeq \F^{G\times H}(P) \simeq \F^H(B)$ 
via pullback $p^*$ and $q^*$ (Lemma 3.3 in \cite{O1}). 
We thus define $A_*(\X):=A_{*+g}^G(X)$ and $\F(\X):=F^G(X)$, 
also $\F_{inv}(\X):=\F^G_{inv}(X)$, through the canonical identification. 

Given proper representable morphisms of quotient stacks $f: \X \to \Y$ 
and any presentations $\X=[X/G]$, $\Y=[Y/H]$, we define 
the pushforward $f_*: A_*(\X) \to A_*(\Y)$ by 
$$f_*^H\circ (q^{*})^{-1}\circ p^*: A_{n+g}^G(X) \to A_{n+h}^H(Y)$$
and also  
$f_*: \F(\X) \to \F(\Y)$ in the same way. 
By the identification $(q^{*})^{-1}\circ p^*$, 
everything is reduced to the equivariant setting (the diagram (\ref{d1})).  

\begin{lemma}  The above $\F$ and $A_*$ satisfy the following properties:    \\
{\rm (i)} For proper representable morphisms of quotient stacks $f$,   
the pushforward $f_*$ is well-defined; \\
{\rm (ii)}  Let $f_1: \X_1 \to \X_2$, $f_2: \X_2 \to \X_3$ and 
$f_3:\X_1 \to \X_3$ be proper representable morphisms of stacks 
so that $f_2\circ f_1$ is isomorphic to $f_3$, 
then  $f_{2*}\circ f_{1*}$  is isomorphic to $f_{3*}$ 
{\rm (}$f_{3*}=f_{2*}\circ f_{1*}$ 
using a notational convention in Remark 5.3, \cite{Gomez}{\rm )}. 

\end{lemma}

\proof
Look at the diagram below, where $\X_i=[X_i/G_i]$ ($i=1,2,3$). 
We may regard $\X_1=[X_1/G_1]=[B_1/G_2]=[B_3/G_3]$, and so on. 
{\rm (i)}  Put $f=f_1$, then the well-definedness of the pushforward $f_{1*}$ (in both of $\F$ and $A_*$) 
is easily checked by taking fiber products and by the canonical identification. 
{\rm (ii)}  Assume that there exists 
an isomorphism of functors $\alpha: f_2\circ f_1 \to f_3$ 
(i.e., a $2$-isomorphism of $1$-morphisms). 
Then two $G_3$-equivariant morphisms 
 $\bar{f}_2\circ \bar{f}_1$ 
and $\bar{f}_3$ from $B_3$ to $X_3$ 
coincide up to isomorphisms of $B_3$ and of $X_3$ 
which are encoded in the definition of $\alpha$, 
hence their $G_3$-pushforwards coincide up to the chosen isomorphisms. 
\qed

\

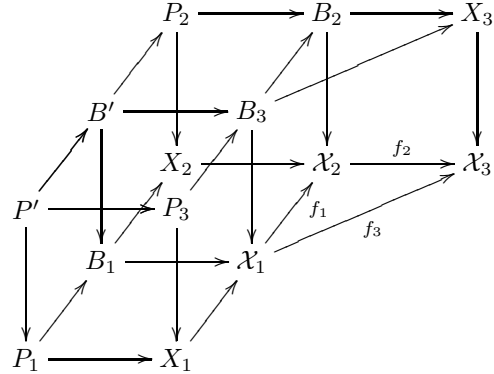
\begin{center}
{\small 
\begin{xy}
(0,0)*+{   }; 
(25,0) *+{P_1}="P1"; (45,0) *+{X_1}="X1";
(35,13) *+{B_1}="B1"; (55,13) *+{\X_1}="XX1"; 
(25,20) *+{P'}="Pd"; (45,20) *+{P_3}="P3"; 
(45,26) *+{X_2}="X2"; (65,26) *+{\X_2}="XX2"; (85,26) *+{\X_3}="XX3";
(35,33) *+{B'}="Bd"; (55,33) *+{B_3}="B3"; 
(45,46) *+{P_2}="P2"; (65,46) *+{B_2}="B2"; (85,46) *+{X_3}="X3";
(64,20) *+{{}_{f_1}};
{\ar@{->} "P1";"X1"}; {\ar@{->} "P1";"B1"};
{\ar@{->} "X1";"XX1"};
{\ar@{->} "Pd";"P1"};{\ar@{->} "Pd";"P3"};{\ar@{->} "Pd";"Bd"};
{\ar@{->} "P3";"X1"};{\ar@{->} "P3";"B3"};
{\ar@{->} "Pd";"P3"};{\ar@{->} "Pd";"Bd"};
{\ar@{->} "B1";"XX1"};{\ar@{->} "B1";"X2"};
{\ar@{->} "Bd";"B3"};{\ar@{->} "Bd";"P2"};{\ar@{->} "Bd";"B1"};
{\ar@{->} "XX1";"XX2"};{\ar@{->}_{f_3} "XX1";"XX3"};
{\ar@{->} "X2";"XX2"};
{\ar@{->}^{f_2} "XX2";"XX3"};
{\ar@{->} "B3";"X3"};{\ar@{->} "B3";"B2"};{\ar@{->} "B3";"XX1"};
{\ar@{->} "P2";"X2"};{\ar@{->} "P2";"B2"};
{\ar@{->} "B2";"X3"};{\ar@{->} "B2";"XX2"};
{\ar@{->} "X3";"XX3"};
\end{xy}
}
\end{center}

\subsection{Chern-MacPherson transformation} 
We assume that $X$ is a quasi-projective scheme or algebraic space 
with action of $G$. 
Then $X_G$ exists as an algebraic space, 
hence $C_*(X_G)$ makes sense. 
Take the vector bundle $TU_G:=X\times_G (U\oplus V)$ over $X_G$, 
 i.e., the pullback 
of the tautological vector bundle $(U\times V)/G$ over $U_G$ 
induced by the projection $X_G \to U_G$. 
Our natural transformation 
$$C^{G}_*: \F^{G}(X) \to A_*^{G}(X)$$ 
is defined to be the inductive limit of  
$$
T_{U, *}:=c(TU_G)^{-1} \frown C_*: \F(X_G) \to A_*(X_G)
$$
over the direct system of representations of $G$, 
see \cite{O1} for the detail. 

Roughly speaking, the $G$-equivariant CSM class 
$C_*^G(X)\, (:=C_*^G(\jeden_X))$ looks like 
``$c(T_{BG})^{-1}\frown C_*(EG\times_G X)$", 
where $EG\times_G X \to BG$ means the universal bundle (as ind-schemes) 
with fiber $X$ and group $G$, 
that has been justified using a different inductive limit of Chow groups, 
 see Remark 3.3 in \cite{O1}. 

\begin{lemma}
{\rm (i)}  In the same notation as in the diagram (\ref{d2}) in 2.1, 
the following diagram commutes: 
 $$\xymatrix{
\F^{G}(X)  \ar[d]_{C^{G}_*} \ar[r]^{p^{*}\;}_{\simeq\;\;} & \F^{G\times H}(P)  \ar[d]^{C^{G\times H}_*}    \\
A^{G}_{*+g}(X) \ar[r]_{p^*\;}^{\simeq\;\;} & A^{G\times H}_{*+g+h}(P)  
}
$$
{\rm (ii)}  In particular, $C_*: F(\X) \to A_*(\X)$ is well-defined. \\
{\rm (iii)}  $C_*f_*=f_*C_*$ for proper representable morphisms $f: \X \to \Y$. 
\end{lemma}

\proof 
{\rm (i)} This is essentially the same as  Lemma 3.1 in \cite{O1} which 
shows the well-definedness of $C_*^G$. 
Apply 
the Verdier-Riemann-Roch \cite{Yokura} to  
the projection of the affine bundle $\bar{p}: P_{G\times H} \to X_G$ (with fiber $U_2$), 
then we have the following commutative diagram 
 $$\xymatrix{
\F(X_G)  \ar[d]_{C_*} \ar[r]^{\bar{p}^{*}\;}_{\;\;} & \F(P_{G\times H})  \ar[d]^{C_*}    \\
A_{*+\ell_1}(X_G) \ar[r]_{\bar{p}^{**}\;}^{\;\;} & A_{*+\ell_1+\ell_2}(P_{G\times H})  
}
$$
where $\bar{p}^{**}=c(T_{\bar{p}}) \frown \bar{p}^*$ 
and $T_{\bar{p}}$ is the relative tangent bundle of $\bar{p}$. 
The twisting factor $c(T_{\bar{p}})$ in $\bar{p}^{**}$ is 
cancelled by the factors in $T_{U_1,*}$ 
and $T_{U,*}$: 
In fact, since $T_{\bar{p}} = \bar{q}^*TU_{2H}$,  $T_{\bar{q}} = \bar{p}^*TU_{1G}$ and 
$$TU_{G\times H}= P \times_{G\times H} (T(U_1\oplus U_2)) 
= T_{\bar{p}}\oplus T_{\bar{q}},$$
we have 
\begin{eqnarray*}
T_{U,*} \circ \bar{p}^* (\alpha)&=& c(TU_{G\times H})^{-1} \frown C_*(\bar{p}^*\alpha)\\
&=&
c(T_{\bar{p}}\oplus T_{\bar{q}})^{-1} c(T_{\bar{p}}) \frown \bar{p}^*C_*(\alpha)\\
&=&
c(T_{\bar{q}})^{-1}  \frown \bar{p}^*C_*(\alpha)\\
&=& 
 \bar{p}^*(c(T{U_1}_G)^{-1}  \frown C_*(\alpha)) \\
 &=&  \bar{p}^* \circ  T_{U_1, *} (\alpha). 
\end{eqnarray*}
Taking the inductive limit, we conclude that $C_*^{G\times H} \circ p^*=p^* \circ C_*^{G}$. 
Thus {\rm (i)}  is proved. 
The claim {\rm (ii)}  follows from {\rm (i)} . 
By {\rm (ii)} , we may consider $C_*$ as the $H$-equivariant 
Chern-MacPherson transformation $C_*^H$ given in  \cite{O1}, 
thus {\rm (iii)}  immediately follows from the naturality of $C_*^H$. 
\qed 

\

The above lemmas show the following theorem 
(cf. Theorem 3.5, \cite{O1}): 

\begin{theorem}\label{theorem} 
Let $\mathcal{C}$ be the category whose objects are 
{\rm (}possibly non-separated{\rm )} Artin quotient stacks $\mathcal{X}$ 
having the form $[X/G]$ of  
separated algebraic spaces $X$ of finite type with action   
of smooth linear algebraic groups $G$; 
morphisms in $\mathcal{C}$ are assumed to be proper and representable. 
Then for the category $\mathcal{C}$, 
we have a unique natural transformation $C_*: \F(\X) \to A_*(\X)$ 
with integer coefficients 
so that it coincides with the ordinary MacPherson transformation 
when restricted to the category of quasi-projective schemes. 
\end{theorem}

\subsection{Degree} 
Let $g=\dim G$.  The $G$-classifying stack $BG=[pt/G]$ has 
(non-positive) virtual dimension $-g$, hence 
$$A_{-n}(BG)=A_{-n+g}^G(pt)=A^{n-g}_G(pt)=A^{n-g}(BG)$$
 for any integer $n$ (trivial for $n<g$). 
 We often use this identification. 
In particular, $A_{-g}(BG)=A^0(BG)=\Z$. 

Let $\X=[X/G]$ in $\mathcal{C}$ with $X$ projective and equidimensional of dimension $n$. 
Then we can take the representable morphism $pt: \X \to BG$: 
$$
\xymatrix{
G\times X \ar[r]^{q} \ar[d]_{p} & X \ar[r]^{\bar{pt}} \ar[d]^{u}  & pt \ar[d] \\
X \ar[r]_{u} & \X \ar[r]_{pt} & BG
}
$$
Here are some remarks: 
\begin{enumerate}
\item[(i)] 
For a $G$-invariant function $\alpha \in \F_{inv}(\X)=\F^G_{inv}(X)$, 
it is obvious that $(q^*)^{-1}\circ p^*(\alpha)=\alpha$, hence we have 
$$pt_*(\alpha) =  \bar{pt}_*(q^*)^{-1} p^*(\alpha)= \bar{pt}_*\alpha=\int_X\alpha =\chi(X; \alpha),$$ 
which is called the {\it integral}, or {\it weighted Euler characteristic} of the invariant function $\alpha$. 
In particular, by the naturality, $pt_*C_*(\alpha)=C_*(pt_*\alpha)=\chi(X; \alpha)$. 
More generally, in \cite{O1} we have defined the $G$-degree of 
{\it equivariant constructible function} $\alpha \in \F(\X)$ 
by $pt_*(\alpha) \in \F^G(pt)=\F(BG)$, which is a 
 `constructible' function over $BG$. Then 
$pt_*C_*(\alpha)=C_*(pt_*\alpha) \in A^*(BG)$, 
being a polynomial or power series in universal $G$-characteristic classes. 

\item[(ii)] 
For invariant functions $\alpha \in \F_{inv}(\X)$ and for $i<-g$ and  $i>n-g$,  
the $i$-th component $C_i(\alpha)$ is trivial. 
A possibly nontrivial highest degree term 
$C_{n-g}(\alpha)  \in A_{n-g}(\X) \, (=A_{n}^G(X))$ is a linear sum of 
 the $G$-fundamental classes $[X_i]_G$ 
of irreducible components $X_i$  (the virtual fundamental class of dimension $n-g$) . 
As a notational convention, 
let $\jeden_\X^{(0)}$ denote the constant function $\jeden_X \in \F^G_{inv}(X) \, =\F_{inv}(\X)$ for a presentation $\X=[X/G]$. In particular, 
if $X$ is smooth, then 
$$C_*(\jeden_\X^{(0)})=C^G_*(\jeden_X)=c^G(TX) \frown [X]_G \in A_{*+g}^G(X)=A_*(\X).$$ 

\item[(iii)]  
From the viewpoint of the enumerative theory in classical projective algebraic geometry 
(e.g. see \cite{P}), a typical type of degrees often arises in the following form:  
$$\int \, pt_*(c(E) \frown C_*(\alpha)) \in A^0(BG)$$
for  some vector bundle $E$ over $\X$ and a constructible function $\alpha \in \F _{inv}(\X)$. 
\end{enumerate}

\section{Deligne-Mumford stacks}
It would be meaningful   
to restrict $C_*$ to a subcategory of certain quotient stacks having finite stabilizer groups, 
which form a reasonable class of Deligne-Mumford stacks  
(including smooth DM stacks). 

\begin{theorem} 
Let $\mathcal{C}_{\rm DM}$ be the category  of 
Deligne-Mumford stacks of finite type which admits a  locally closed embedding into 
some smooth proper DM stack with projective coarse moduli space: morphisms in $\mathcal{C}_{\rm DM}$ are assumed to be proper and representable. 
 Then for $\mathcal{C}_{\rm DM}$ 
 there is a unique natural transformation $C_*: \F(\X) \to A_*(\X)$ 
 satisfying the normalization property: $C_*(\jeden_X)=c(TX)\frown [X]$ for  smooth schemes. 
 \end{theorem}
 
This is due to Theorem 5.3 in Kresch \cite{Kresch} which states that 
a DM stack in $\mathcal{C}_{\rm DM}$ is in fact realized by a quotient stack  in $\mathcal C$. 
In \cite{Kresch}, such a DM stack is called to be ({\it quasi-}){\it projective}. 

\begin{remark}
{\rm 
(i) In the above theorem,  
the embeddability into smooth stack 
(or equivalently the resolution property in \cite{Kresch}) is required, 
that seems natural, since original MacPherson's theorem 
requires such a condition \cite{Mac, Ken}. 
In order to extend $C_*$ for more general Artin stacks with values in Kresch's Chow groups,  
we need to find some technical gluing property. \\
(ii) 
We may admit proper {\it non-representable} morphisms of DM stacks 
if we use rational coefficients.  In fact for such morphisms  
the pushforward of Chow groups with rational coefficients is defined \cite{Vistoli}. 
}
\end{remark}

\subsection{Modified pushforwards} \label{mod}

The theory of constructible functions for Artin stacks has been established by Joyce \cite{Joyce}. 
Below let us work with $\Q$-valued constructible functions and Chow groups with $\Q$-coefficients. 
For stacks $\X$ in $\mathcal{C}_{\rm DM}$, 
each geometric point  $x: pt= {\rm Spec}\, k \to \X$ has a finite stabilizer group $Aut(x)(={\rm Iso}_x(x,x))$. 
Then the group of  constructible functions  $\underline{\alpha}$ in the sense \cite{Joyce} 
is canonically identified with 
the subgroup $F_{inv}(\X)_\Q=F^G_{inv}(X)_\Q$  
of invariant constructible functions $\alpha$ over $X$ 
in the following way (the bar indicates constructible functions over the set of all geometric points  $\X(k)$): 
For each $k$-point $x: pt \to \X$, the value of $\alpha$ over the orbit $x\times_\X X$ is given by 
$|Aut(x)|\cdot \underline{\alpha}(x)$, that is, 
$$\F(\X(k))_\Q \simeq \F_{inv}(\X)_\Q  \;\; (\; \subset \F(\X)_\Q \;), \quad \underline{\alpha}  \leftrightarrow 
\alpha=\jeden_\X \cdot \pi^*\underline{\alpha}, $$
where $\pi$ is the projection to $\X(k)$, 
 $\alpha\cdot \beta$ is  the canonical multiplication on $\F(\X)_\Q$, 
 $(\alpha\cdot \beta)(x):=\alpha(x)\beta(x)$, and 
 $$\jeden_\X:=|Aut(\pi(-))| \in \F_{inv}(\X)_\Q.$$

It is shown by Tseng \cite{Tseng} that 
if $\X$ is a smooth DM stack, 
$C_*(\jeden_{\X})$ coincides with (pushforward of the dual to) the total Chern class 
 of the tangent bundle of the corresponding smooth inertia stack. 

From a viewpoint of classical group theory, 
it would be natural to measure {\it how large of the stabilizer group is} 
by comparing it with a fixed group $A$, 
that leads us to define a $\Q$-valued constructible function over $\X(k)$. 
Here the group $A$ is supposed to be, e.g., a finitely generated Abelian group 
(we basically consider $A=\Z^m$, $\Z/r\Z$, etc). 
Accordingly to \cite{O1, O2}, we define 
{\it the canonical constructible function measured by group $A$} 
which assigns to any geometric point $x$ 
the number of group homomorphisms of $A$ into $Aut(x)$: 
$$\underline{\jeden_{\X}^{A}}(x):=\frac{|\, {\rm Hom}\, (A, Aut(x))\, |}{|\, Aut(x)\, |}  \; \in \Q.$$
The corresponding invariant constructible function is denoted by 
$\jeden_\X^A \in \F_{inv}(\X)_\Q$, or  often by 
$\jeden_{X;G}^{A} \in \F^G_{inv}(X)_\Q$  
when a presentation $\X=[X/G]$ is specified. Namely, 
the value of $\jeden_{X;G}^{A}$ on the $G$-orbit expressed by 
$x: pt \to \X$ is $|\, {\rm Hom}\, (A, Aut(x))\, |$. 
The function for $A=\Z$ is nothing but $\jeden_\X$ 
in our convention, and for $A=\{0\}$ it is $\jeden_\X^{(0)}=1$. 
If $A=\Z^2$,  the function counts the number of mutually commuting pairs in $Aut(x)$, 
hence its integral corresponds to 
the orbiforld Euler number (in physicist's sense), see \cite{O2}.

Define  
$T^A_\X: \F(\X)_\Q \to \F(\X)_\Q$ by 
the multiplication  $T^A_\X(\alpha):=\jeden_{X;G}^{A} \cdot \alpha$. 
This is a $\Q$-algebra isomorphism, 
for $\jeden_{X;G}^{A}$ is an unit in $\F(\X)_\Q$. 
A new pushforward is introduced 
for proper representable morphisms $f:\X \to \Y$ in $\mathcal{C}_{\rm DM}$ by
$$f^{A}_*: \F(\X)_\Q \to \F(\Y)_\Q,  \quad 
\alpha \mapsto (T^A_\Y)^{-1} \circ f_*\circ T^A_\X(\alpha).$$
Obviously, 
  $g^{A}_* \circ f^{A}_*=(g\circ f)^{A}_*$. 
 The following theorem says that there are 
 several variations of theories of integration with values in Chow groups 
for Deligne-Mumford stacks:

\begin{theorem} Given a finitely generated Abelian group $A$, 
let $F^{A}$ denote the new covariant functor of constructible functions 
for the category  $\mathcal{C}_{\rm DM}$, given by 
$F^{A}(\X)_\Q:=F(\X)_\Q$ and  the pushforward  by  $f^{A}_*$. 
Then, 
 $C^{A}_*:=C_*\circ T^A_\X: F^{A}(\X)_\Q \to A_*(\X)_\Q$ is a natural transformation. 
 \end{theorem}
 
\proof 
It is straightforward that 
 $f_*\circ C^{A}_*=f_*\circ C_*\circ T^A_\X = C_* \circ f_*\circ T^A_\X
 = C_* \circ T^A_\Y \circ (T^A_\Y)^{-1} \circ f_*\circ T^A_\X
=C_*^{A} \circ f_*^{A}.$ \qed

\section{Other characteristic classes}
The method in the preceeding sections
is applicable to other characteristic classes
(over $\C$ or a field $k$ of characteristic $0$). 

As the most general additive characteristic class for
singular varieties, {\it the Hirzebruch class transformation} 
$$T_{y*}: K_0(Var/X) \to A_*(X) \otimes \Q[y]$$
was recently introduced by Brasselet-Sch\"urmann-Yokura \cite{BSY}: 
For possibly singular varieties $X$ (and proper morphisms between them), 
$T_{y*}$ is a unique natural transformation from 
the Grothendieck group $K_0(Var/X)$ of the monoid of 
isomorphism classes of morphisms $V \to X$ to 
the rational Chow group of $X$ with a parameter $y$ such that it satisfies that 
$$T_{y*}[X\stackrel{id}{\to}X]=\widetilde{td_y}(TX) \frown [X], \quad 
\mbox{for smooth $X$,}$$
where $\widetilde{td_y}(E)$ denotes 
the modified Todd class of vector bundles: 
$$\widetilde{td_y}(E)=\prod_{i=1}^r \left(\frac{a_i(1+y)}{1-e^{-a_i(1+y)}}-a_iy \right),$$
when $c(E)=\prod_{i=1}^r (1+a_i)$, see  \cite{BSY, SY}. 
Note that the associated genus is well-known Hirzebruch's $\chi_y$-genus, 
which specializes to: 
the Euler characteristic if $y=-1$,  
the arithmetic genus if $y=0$, 
and the signature if $y=1$. 
Hence, $T_{y*}$ gives 
a generalization of the $\chi_y$-genus to 
homology characteristic class of singular varieties, 
which unifies the following 
singular Riemann-Roch type formulas in canonical ways: 
\begin{itemize}
\item 
($y=-1$) 
the Chern-MacPherson transformation $C_*$ \cite{Mac, Ken}; 
\item 
($y=0$)   
Baum-Fulton-MacPherson's Todd class transformation $\tau$ \cite{BFM}; 
\item 
($y=1$)   
Cappell-Shaneson's homology $L$-class transformation $L_*$ \cite{CS}. 
\end{itemize}

For a quotient stack $\X=[X/G] \in \mathcal{C}$ in Theorem 1,
we denote by $K_0(\mathcal{C}/\X)$
the Grothendieck group of the monoid
of isomorphism classes of representable morphisms
of quotient stacks to the stack $\X$. 
To each element $[\mathcal{V} \to \X] \in K_0(\mathcal{C}/\X)$,
we take 
a $G$-equivariant morphism $V \to X$ 
where $V:=\mathcal{V}\times_{\X} X$ with natural $G$-action 
so that $\mathcal{V}=[V/G]$, and associate 
 a class of morphisms of algebraic spaces 
$[V_G \to X_G] \in K_0(Var/X_G)$.
We then define
$$T_{y*}: K_0(\mathcal{C}/\X) \to A_*(\X) \otimes \Q[y]$$
by assigning to $[\mathcal{V} \to \X]$ the inductive limit 
(over all $G$-representations) of 
$$\widetilde{td_y}^{-1}(TU_G) \frown T_{y*}[V_G \to X_G] \in 
A_*(X_G)\otimes \Q[y].$$
This is well-defined, because the Verdier-Riemann-Roch for $T_{y*}$ holds 
(Corollary 3.1 in \cite{BSY}) 
and the same proof of Lemma 2 can be used in this setting. 
Note that in each degree of grading 
the limit stabilizes, thus the coefficient is a polynomial in $y$. 
So we obtain an extension of $T_{y*}$ to the category $\mathcal{C}$,  
and hence also to $\mathcal{C}_{\rm DM}$. 

It turns out  that  at special values $y=0, \pm 1$, 
$T_{y*}$ corresponds to: 
\begin{itemize}
\item 
($y=-1$)  the $G$-equivariant Chern-MacPherson transformation \cite{O1},  
i.e., $C_*$ as described in section 2 above; 
\item 
($y=0$)  
the $G$-equivariant Todd class transformation \cite{EG, BZ}, 
given by the limit of $td^{-1}(TU_G) \frown \tau$; 
\item 
($y=1$)  the $G$-equivariant singular $L$-class transformation 
given by the limit of
$(L^*)^{-1}(TU_G) \frown L_*$, where $L^*$ is 
the (cohomology) Hirzebruch-Thom $L$-class. 
\end{itemize}
Applications will be considered in another paper.

\end{document}